\documentclass[12pt]{amsart}
\usepackage{amscd,amssymb}
\usepackage[graph,frame,poly,arc]{xy}  
\usepackage[plainpages,backref,urlcolor=blue]{hyperref}

\theoremstyle{plain}
\newtheorem{thm}[subsection]{Theorem}

\newtheorem{prop}[subsection]{Proposition}
\newtheorem{cor}[subsection]{Corollary}

\theoremstyle{definition}
\newtheorem{rk}[subsection]{Remark}

\newtheorem{ex}[subsection]{Example}
\newtheorem{conj}[subsection]{Conjecture}

\numberwithin{equation}{section}
\newcommand{\A}{{\mathcal A}}

\newcommand{\R}{\mathbb{R}}
\newcommand{\C}{\mathbb{C}}
\newcommand{\PP}{\mathbb{P}}

\newcommand{\N}{\mathbb{N}}

\begin{document}
%\date{June 4, 2009}

\title [Rational cuspidal  curves and local cohomology of Jacobian rings]
{On rational cuspidal plane curves, and the local cohomology of Jacobian rings}

\author[Alexandru Dimca]{Alexandru Dimca$^1$}
\address{Universit\'e C\^ ote d'Azur, CNRS, LJAD, France }
\email{dimca@unice.fr}

%\thanks{$^1$ Partially supported by Institut Universitaire de France.} 

\thanks{$^1$ This work has been supported by the French government, through the $\rm UCA^{\rm JEDI}$ Investments in the Future project managed by the National Research Agency (ANR) with the reference number ANR-15-IDEX-01.}



\begin{abstract} This note gives the complete projective classification of rational, cuspidal plane curves of degree at least 6, and having only weighted homogeneous singularities. 
It also sheds new light on some previous characterizations of free and nearly free curves in terms of Tjurina numbers. Finally, we suggest a stronger form of Terao's conjecture on the freeness of a line arrangement being determined by its combinatorics.

\end{abstract}
 
\maketitle

%\tableofcontents

\section{Introduction} 

The  main result we prove in this note is the following.

\begin{thm}
\label{thmRCC}
 Let $C$ be a rational, cuspidal plane curve of degree $d \geq 6$ such that $C$ has only weighted homogeneous singularities. Then $C$ is projectively equivalent to exactly one of the following $\phi(d)/2$ models
$$C_{d,k}: y^d+x^kz^{d-k}=0,$$
where the integer $k$ satisfies $1 \leq k <d/2$ and $k$ is relatively prime to $d$. Here $\phi$ denotes the Euler function, with $\phi(d)$ counting the number of integers $m$, where $1 \leq m \leq d-1$ and $m$ relatively prime to $d$.
\end{thm}

A slightly stronger result is stated and proved below, see Propositions \ref{propRCC} and  \ref{propRCC2}.
For detailed information on the rational, cuspidal plane curves of degree $\leq 5$, we refer to \cite{Moe} and the references given there. The curves $C_{d,k}$ above are called {\it binomial cuspidal curves}  in \cite[Section 7.1]{Moe}. The necessary restriction $d \geq 6$ in Theorem \ref{thmRCC} is discussed in Example \ref{exMoe} below.

Theorem \ref{thmRCC} is quite surprising, given that the classification of rational, cuspidal plane curves, even of those having at most two singularities, is rather complicated, see Sakai-Tono paper \cite{SaTo}, or Propositions
3.1 and 3.2 in \cite{DStFD} where their results are quoted.
It is known that there are no rational, cuspidal {\it free} curves of degree $d \geq 6$ having only weighted homogeneous singularities, see \cite[Theorem 2.8]{DStFD}. Theorem \ref{thmRCC} says that if the word {\it free} is deleted from the above statement, only the binomial cuspidal curves may occur, and they are all nearly free. Note also that a weighted homogeneous cusp has only one Puiseux pair, but the class of plane curve singularities with exactly one Puiseux pair is much larger. This explains the complicated classification of rational unicuspidal curves with a unique Puiseux pair in \cite{FLMN}.

To prove Theorem \ref{thmRCC}, we use a number of key results, proved by A. A. du Plessis and C.T.C. Wall in \cite{duPCTC, duPCTC2}, the main one being stated below in Theorem \ref{thmCTC}. To make our note more self-contained, we include a short proof of a similar result to Theorem \ref{thmCTC},
which we explain now.

Let $S=\C[x,y,z]$ be the graded polynomial ring in three variables $x,y,z$ and let $C:f=0$ be a reduced curve of degree $d$ in the complex projective plane $\PP^2$. The minimal degree of a Jacobian relation, or Jacobian syzygy, for the polynomial $f$ is the integer $mdr(f)$
defined to be the smallest integer $m\geq 0$ such that there is a nontrivial relation
\begin{equation}
\label{rel_m}
 af_x+bf_y+cf_z=0
\end{equation}
among the partial derivatives $f_x, f_y$ and $f_z$ of $f$ with coefficients $a,b,c$ in $S_m$, the vector space of  homogeneous polynomials of degree $m$. When $mdr(f)=0$, then $C$ is a union of lines passing through one point, a situation easy to analyse. We assume from now on that 
$$ mdr(f)\geq 1.$$
Denote by $\tau(C)$ the global Tjurina number of the curve $C$, which is the sum of the Tjurina numbers of the singular points of $C$.
We denote by $J_f$ the Jacobian ideal of $f$, i.e. the homogeneous ideal in $S$ spanned by $f_x,f_y,f_z$, and  by $M(f)=S/J_f$ the corresponding graded ring, called the Jacobian (or Milnor) algebra of $f$.
 Let $I_f$ denote the saturation of the ideal $J_f$ with respect to the maximal ideal ${\bf m}=(x,y,z)$ in $S$ and consider the local cohomology group 
 $$N(f)=I_f/J_f=H^0_{\bf m}(M(f).$$
 
 It was shown in \cite{DPop} that the graded $S$-module  $N(f)$ satisfies a Lefschetz type property with respect to multiplication by generic linear forms. This implies in particular the inequalities
$$0 \leq n(f)_0 \leq n(f)_1 \leq ...\leq n(f)_{[T/2]} \geq n(f)_{[T/2]+1} \geq ...\geq n(f)_T \geq 0,$$
where $T=3d-6$ and $n(f)_k=\dim N(f)_k$ for any integer $k$. We set
$$\nu(C)=\max _j \{n(f)_j\}.$$
When $d=2m$ is even, then the above implies that
$n(f)_{3m-3}=\nu(C).$
When $d=2m+1$ is odd, then the above and the self duality of the graded $S$-module $N(f)$, see \cite{Se, SW}, implies that
$$n(f)_{3m-2}=n(f)_{3m-1}=\nu(C).$$

The second main result of this note is the following.

\begin{thm}
\label{thmN}
Let $C:f=0$ be a reduced plane curve of degree $d$ and let $r=mdr(f)$.
Then the following hold.
\begin{enumerate}
\item If $r < d/2$, then
$$\nu(C)=(d-1)^2-r(d-1-r)-\tau(C).$$

\item If $r \geq (d-2)/2$, then 
$$\nu(C)= \lceil   \  \frac{3}{4}(d-1)^2 \    \rceil  -\tau (C).$$

\end{enumerate}
\end{thm}
Here, for any real number $u$,  $\lceil     u    \rceil $ denotes the round up of $u$, namely the smallest integer $U$ such that $U \geq u$.
Written down explicitly, this means that for $d=2m$ is even and $r \geq m-1$, one has
$\nu(C)=3m^2-3m+1-\tau(C)$, while for $d=2m+1$ is odd and  $r \geq m$, one has
$\nu(C)=3m^2-\tau(C).$ For $(d-2)/2 \leq r <d/2$, both formulas in (1) and (2) apply, and they give the same result for $\nu(C)$.
The relation between Theorem \ref{thmN} and du Plessis-Wall result in Theorem \ref{thmCTC} is discussed in Remark \ref{rkCTC} below.

When $C:f=0$ is a line arrangement, examples due to G. Ziegler show that 
the invariant $mdr(f)$ is not combinatorially determined, see \cite[Example 4.3]{AD}. The above result suggests that the following stronger version of H. Terao's conjecture, saying that the freeness of a line arrangement is combinatorially determined, might be true. For more on Terao's conjecture and free hyperplane arrangements we refer to \cite{DHA}.

\begin{conj}
\label{conj1}
Let $C:f=0$ be a line arrangement in $\PP^2$. Then the invariant $\nu(C)$ is combinatorially determined.
\end{conj}

For a reduced plane curve one may state the following.
\begin{conj}
\label{conj1G}
Let $C:f=0$ be a reduced plane curve $\PP^2$. Then the invariant $\nu(C)$ is determined by the degree of $C$ and the list of the analytic types of the isolated singularities of $C$.
\end{conj}

Some cases where these conjecture hold are described in Example \ref{exMoe} (i) (where the corresponding invariant $mdr(f)$ is not determined by the data stated in Conjecture  \ref{conj1G}), and in 
Proposition \ref{propTeraoC}.

\section{The proof of Theorem \ref{thmN}}

 Consider the graded $S-$submodule $AR(f) \subset S^{3}$ of {\it all relations} involving the partial derivatives of $f$, namely
$$\rho=(a,b,c) \in AR(f)_m$$
if and only if  $af_x+bf_y+cf_z=0$ and $a,b,c$ are in $S_m$. We set $ar(f)_k=\dim AR(f)_k$ for any integer $k$.

Consider the rank two vector bundle $T\langle C\rangle=Der(-{\rm log} C)$ of logarithmic vector fields  along $C$, which is the coherent sheaf associated to the graded $S$-module $AR(f)(1)$.
Using the results in the third section of \cite{DS14}, for any integer $k$ one has
\begin{equation}
\label{chi}
\chi(T\langle C\rangle(k))=3{k+3 \choose 2}-{d+k+2 \choose 2} +\tau(C).
\end{equation}
Moreover, one has the following for $E=T\langle C\rangle$ and any integer $k$, see \cite{DS14}, \cite{Se}.
\begin{equation}
\label{chi2}
 h^0((E(k))=ar(f)_{k+1}, \  \  h^1((E(k))=n(f)_{d+k} \text{ and
}  h^2((E(k))=ar(f)_{d-5-k}.
\end{equation}

Assume that we are in the case $d=2m$ and
apply the formulas  \eqref{chi} and \eqref{chi2} for  $k=m-3$. We get
$$2ar(f)_{m-2}-\nu(C)=\tau(C)-3m^2+3m-1.$$
Let $r=mdr(f)$ and note that $ar(f)_{m-2}=0$ if $r \geq m-1=[(d-1)/2]$
and 
$$ar(f)_{m-2}= \dim S_{m-2-r}={m-r \choose 2}$$
if $r \leq m-2$. Indeed, it follows from \cite[Lemma 1.4]{ST} that the $S$-module $AR(f)$ cannot have two independent elemens of degree $d_1$ and $d_2$ satisfying $d_1+d_2 <d-1$.
Assume next that we are in the case $d=2m+1$ and
apply the formulas  \eqref{chi} and \eqref{chi2} again for  $k=m-3$. We get
$$ar(f)_{m-2}-\nu(C)+ar(f)_{m-1}=\tau(C)-3m^2.$$
As above, if  $r=mdr(f)$, we have $ar(f)_{m-2}=ar(f)_{m-1}=0$ if $r \geq m=[(d-1)/2]$
and 
$$ar(f)_{m-2}= \dim S_{m-2-r}={m-r \choose 2} \text{ and } ar(f)_{m-1}= \dim S_{m-1-r}={m-r +1\choose 2} $$
if $r \leq m-1$, with the convention ${1 \choose 2}=0$. 
These formulas, plus some simple computation, prove Theorem \ref{thmN} stated in the Introduction.

\begin{rk}
\label{rkthmN}
An alternative proof of Theorem \ref{thmN} can be obtained using the results in \cite{AD}. Indeed, Proposition 3.2 in \cite{AD} implies that 
$d_1^{L_0}=r$ for $r <(d-2)/2$ and $d_1^{L_0}=\lfloor (d-1)/2 \rfloor $ for $r  \geq (d-2)/2$. Here $(d_1^{L_0},d_2^{L_0})$, with  $d_1^{L_0} \leq d_2^{L_0}$, denotes the splitting type of the vector bundle $T\langle C\rangle(-1)$ along a generic line $L_0$ in $\PP^2$. Then \cite[Theorem 1.1]{AD} says that 
$$(d-1)^2 -d_1^{L_0}d_2^{L_0}=\tau(C)+ \nu(C).$$
Since $d_1^{L_0}+d_2^{L_0}=d-1$ by \cite[Proposition 3.1]{AD}, these facts give a new proof of Theorem \ref{thmN}. 
\end{rk}
In view of this remark, Conjecture \ref{conj1} may be restated as follows.

\begin{conj}
\label{conj2}
Let $C:f=0$ be a line arrangement in $\PP^2$. Then the generic splitting type  $(d_1^{L_0},d_2^{L_0})$ of the vector bundle $T\langle C\rangle(-1)$ is combinatorially determined.
\end{conj}

\section{Some related results and direct applications} 

We start this section by recalling the following result due to  du Plessis and Wall, see \cite[Theorem 3.2]{duPCTC}.
\begin{thm}
\label{thmCTC}
For positive integers $d$ and $r$, define two new integers by 
$$\tau(d,r)_{min}=(d-1)(d-r-1)  \text{ and } 
\tau(d,r)_{max}= (d-1)^2-r(d-r-1).$$ 
Then, if $C:f=0$ is a reduced curve of degree $d$ in $\PP^2$ and  $r=mdr(f)$,  one has
$$\tau(d,r)_{min} \leq \tau(C) \leq \tau(d,r)_{max}.$$
Moreover, for $r=mdr(f) \geq d/2$, the stronger inequality
$$\tau(C) \leq \tau(d,r)_{max} - {2r+2-d \choose 2}$$
holds.

\end{thm}

\begin{rk}
\label{rkCTC} Let $C:f=0$ be a reduced curve of degree $d$ in $\PP^2$ and  $r=mdr(f)$. Note that the function $\tau(d,r)_{max}$, regarded as a function of $r$, occurs also in Theorem \ref{thmN} (1), which can be restated as
\begin{equation}
\label{eqN}
\tau(C)+\nu(C)=\tau(d,r)_{max},
\end{equation}
for $r=mdr(f) <d/2$.
The inequality
$\tau(C) \leq \tau(d,r)_{max}$ in Theorem \ref{thmCTC} is made more precise, when $r <d/2$, by the result in Theorem \ref{thmN} (1).

On the other hand, note that $\tau(d,r)_{max}$ is a decreasing (resp. increasing) function of $r$ on the interval $[0, (d-1)/2]$ (resp. on the interval $[(d-1)/2, + \infty)$), and the number $ \lceil   \  \frac{3}{4}(d-1)^2 \    \rceil$ that occurs  in Theorem \ref{thmN} (2) is precisely the minimum
of the values of the function $\tau(d,r)_{max}$ for 
$r \in [0, (d-1)/2]\cap \N$. Hence Theorem \ref{thmN} (2) implies only
$$ \tau(C) \leq \lceil   \  \frac{3}{4}(d-1)^2 \    \rceil,$$
for $r \geq d/2$, which is weaker than the inequality in Theorem \ref{thmCTC}, though enough for many applications.

\end{rk}

At the end of the proof of Theorem \ref{thmCTC}, in \cite{duPCTC}, the authors state the following very interesting consequence (of the proof, not of the statement) of Theorem \ref{thmCTC}.
\begin{cor}
\label{corCTC} Let $C:f=0$ be a reduced curve of degree $d$ in $\PP^2$ and  $r=mdr(f)$. One has
$$ \tau(C) =\tau(d,r)_{max}$$
if and only if $C:f=0$ is a free curve, and then $r <d/2$.
\end{cor}
Since a plane curve $C$ is free if and only if $\nu(C)=0$, this characterization of free curves follows also from Theorem \ref{thmN}, as explained in Remark \ref{rkCTC}.

In the paper \cite{Dmax}, we have given an alternative proof of Corollary \ref{corCTC} and have
shown that  a plane curve $C$ is nearly free, which can be defined by the property $\nu(C)=1$, if and only if a similar property holds. Namely, one has the following result, an obvious consequence of Theorem \ref{thmN}.

\begin{prop}
\label{propCTC2} Let $C:f=0$ be a reduced curve of degree $d$ in $\PP^2$ and  $r=mdr(f)$.
One has
$$ \tau(C) =\tau(d,r)_{max}-1$$
if and only if $C:f=0$ is a nearly free curve, and then $r  \leq d/2$.
\end{prop}

Concerning Conjecture \ref{conj1} in Introduction we have the following.

\begin{prop}
\label{propTeraoC} Let $C$ be a plane curve of degree $d$ such that, for $d=2m$ (resp. $d=2m+1$) one has 
$\tau(C) < \tau(d,m-2)_{min}=(m+1)(2m-1)$ (resp. $\tau(C) < \tau(d,m-1)_{min}=2m(m+1)$). Then 
$$\nu(C)= \lceil   \  \frac{3}{4}(d-1)^2 \    \rceil  -\tau (C).$$
In particular, if $C$ be a line arrangement (resp. a reduced plane curve) satisfying these inequalities, then  Conjecture \ref{conj1} (resp. Conjecture \ref{conj1G}) holds for $C$.
\end{prop}

\proof We give the proof only in the case $d=2m$ even. The function
$\tau(d,r)_{min}$ is a strictly decreasing function of $r$ on $\R$, so Theorem \ref{thmCTC} implies that $r=mdr(f) \geq m-1$. But then Theorem \ref{thmN} (2) implies that $\nu(C)$ is determined by $d$ and $\tau(C)$. When $C$ is a line arrangement,  $d$ and $\tau(C)$ are both combinatorially determined. For a plane curve in general, $\tau(C)$ is determined by the analytic type of the local singularities, but not by their topological type.

\endproof

\begin{ex}
\label{exTeraoC} 
Let $C=\A: f=0$ be a line arrangement consisting of $d$ lines. For $d=6$, there are 10 possibilities for the intersection lattice $L(\A)$, see for details \cite[Section 5.6]{DIM}. Out of them, 8 satisfy the inequality $\tau(\A) < 20$ from Proposition \ref{propTeraoC}.
For $d=6$, $r=mdr(f)$ is determined by the intersection lattice $L(\A)$, and there are exactly 8 lattices with $r\geq 2$.

For $d=7$, there are 23 possibilities for the intersection lattice $L(\A)$. Out of them, only 4 satisfy the inequality $\tau(\A) < 24$ from Proposition \ref{propTeraoC}.
For $d=7$, $r=mdr(f)$ is again determined by the intersection lattice $L(\A)$, and there are exactly 19 lattices with $r\geq 3$, for which $\nu(C)=\nu(\A)$ is given by the formula in Proposition \ref{propTeraoC}.

It would be interesting to find a lower bound for $\tau(\A)$ in terms of $d$ and $r$, in the case of line arrangements, which is better than the bound $\tau(r)_{min}$ given by Theorem \ref{thmCTC}.

\end{ex}

\section{The proof of Theorem \ref{thmRCC}} 

The following result is related to the main conjecture in \cite{DStNF}, namely that a rational, cuspidal plane curve is either free or nearly free.  
It is also the first step in proving Theorem \ref{thmRCC}.

Recall  that a plane curve $C$ has only weighted homogeneous singularities if and only if 
$\mu(C)=\tau(C)$, where $\mu(C)$ denotes the sum of all the Milnor numbers of the singularities of $C$, see \cite{KS}. In general one has the obvious inequality
$\tau(C) \leq \mu(C)$.

\begin{prop}
\label{propRCC} Let $C:f=0$ be an irreducible curve in $\PP^2$, of degree $d \geq 6$. Then the following properties are equivalent.
\begin{enumerate}

\item $\tau(C) \geq d^2-4d+8.$

\item $r=mdr(f)=1$

\item  $C$ is a rational, cuspidal plane curve such that  $\mu(C)=\tau(C)$, i.e. $C$ has only weighted homogeneous singularities.
\end{enumerate}
If any of these properties hold, then
$C$ is a nearly free curve.

\end{prop}

\proof To prove that (1) implies (2), recall that the function $\tau(d,r)_{max}$ is a decreasing function of
$r$ on the interval $[1,(d-1)/2]$ and note that $\tau(d,r)_{max}=\tau(d,r')_{max}$ for $r+r'=d-1$. For $r=2$ we get 
$\tau(d,2)_{max}=d^2-4d+7$. Then (1) implies that $r=mdr(f)$ has to be 1. Indeed, the value $r=d-2$ is excluded using the stronger final inequality in Theorem \ref{thmCTC}.
Note that Theorem \ref{thmN} (2) is also enough for this purpose, as discussed in Remark \ref{rkCTC}.
To prove that (2) implies (3), note that
$$\tau(d,1)_{min}=d^2-3d+2 \leq \tau(C) \leq \tau(d,1)_{max}=d^2-3d+3.$$
The genus $g(C)$ is given by
$$2g(C)=(d-1)(d-2) -\sum_p(\mu(C,p)+r(C,p)-1),$$
where $p$ runs through the singular points of $C$, $\mu(C,p)$ is the Milnor number of the singularity $(C,p)$ and $r(C,p)$ is the number of its branches. It follows that
$$\mu(C) =\sum_p\mu(C,p) \leq (d-1)(d-2),$$
with equality if and only if $C$ is a rational cuspidal curve.
Then the  inequality $\tau(C) \leq \mu(C) \leq (d-1)(d-2)$ forces the equality
$$\tau(d,1)_{min}=d^2-3d+2 = \tau(C)$$
and we see that $C$ is nearly free using Proposition \ref{propCTC2}.
The fact that (3) implies (1) is obvious as soon as $d \geq 6$.

\endproof

\begin{rk}
\label{rkTAU}
Irreducible, cuspidal nearly free curves satisfying $\mu(C)=\tau(C)$, but  not rational have been constructed in \cite{B+}. For instance, for any odd integer $k \geq 1$, it is shown that the irreducible curve
$$C_{2k}: f=x^{2k}+y^{2k}+z^{2k}-2(x^ky^k+x^kz^k+y^kz^k)=0,$$
has $3k$ singular points of type $A_{k-1}$ as singularities, it is a nearly free curve with $mdr(f)=k$ and has genus
$$g(C_{2k})=\frac{(k-1)(k-2)}{2}.$$
For $k=3$, we have in this case $\tau(C)=18 <\tau(d,2)_{max}+1=20.$
\end{rk}

The following result completes the proof of Theorem \ref{thmRCC}.

\begin{prop}
\label{propRCC2} Let $C:f=0$ be an irreducible plane curve of degree $d \geq 6$ such that $mdr(f)=1$.
Then $C$ is projectively equivalent to exactly one of the following $\phi(d)/2$ models
$$C_{d,k}: y^d+x^kz^{d-k}=0,$$
where the integer $k$ satisfies $1 \leq k <d/2$ and $k$ is relatively prime to $d$.

\end{prop}

\proof  Using Proposition \ref{propRCC} and \cite[Proposition 1.1]{duPCTC2}, we see that $C$ has a 1-dimensional symmetry, i.e. it admits a 1-dimensional algebraic subgroup $H'$ of $PGL_2(\C)$ as automorphism group. This group lifts uniquely to an algebraic 1-parameter subgroup $H$ of $GL_3(\C)$, not contained in the center and preserving not only the curve $C$, but also its defining equation $f$. Such a subgroup $H$ may be either semi-simple or nilpotent, and
\cite[Proposition 3.1]{duPCTC2} tells us that only the semi-simple case can occur for our situation.
Indeed, as in the proof above, we know that
 $$\tau(C)=(d-1)(d-2)=d^2-3d+2.$$
Note that the statement of  \cite[Proposition 3.1]{duPCTC2} should be slightly corrected, namely the part '' if $y^d$ has non-zero coefficient in $f$ '' is to be replaced by '' if $y^d$ has zero coefficient in $f$ ''.
Once we know that $H$ is semi-simple, we can assume that $H=\C^*$ acting on $\PP^2$ via 
$$t\cdot [x:y:z]=[t^{w_1}x:t^{w_2}y:t^{w_3}z],$$
for some integers $w_j$. Then using the discussion and the notation on page 120 in 
\cite{duPCTC2}, we see that the only possibilities to get irreducible curves correspond to line segments $A\alpha_k$ for $1 \leq k <d/2$ and $k$ relatively prime to $d$. 

\endproof

\begin{ex}
\label{exMoe} 
(i)  Consider the following two rational cuspidal quartics
$$C:f=y^4-xz^3=0 \text{ and } C':f'=y^4-xz^3-y^3z.$$
Then both curves have an $E_6$-singularity located at $[1:0:0]$,
hence they have only weighted homogeneous singularities. However, as noticed in \cite[Section 7.1]{Moe}, in the final part on semi-binomial curves, $C$ and $C'$ are not projectively equivalent. This follows there from the fact that $C$ (resp. $C'$) has one (resp. two) inflection points.
From our point of view, the difference between $C$ and $C'$ is that
$mdr(f)=1$, while $mdr(f')=2$, as a direct computation shows. Note that both curves $C$ and $C'$ are nearly free, since one has
$$\tau(C)=\tau(C')=6=\tau(4,1)_{max}-1=\tau(4,2)_{max}-1.$$
Similar semi-binomial rational cuspidal curves $C:f=0$ of arbitrary degree $d$, with a unique Puiseux pair and $\tau(C) <\mu(C)$, and with all the possible values for $r=mdr(f) \in [2,d/2]$, are discussed in \cite{DStExpo}.

\medskip

\noindent (ii) For $d=5$, any rational cuspidal curve with 4 cusps is projectively equivalent to
$$C:f=16x^4y+128x^2y^2z-4x^3z^2+256y^3z^2-144xyz^3+27z^5=0,$$
see for instance the discussion of the curve $C_8=[(2_3),(2),(2),(2)]$ in
 \cite[Section 6.3]{Moe}.
This curve has three $A_2$ cusps and one $A_6$ cusp, hence only 
weighted homogeneous singularities. A direct computation shows that
$mdr(f)=2$ and that $\tau(C)=\tau(5,2)_{max}$, and therefore $C$ is a free curve.

It is surprizing that the condition $d \geq 6$ in our results above avoids complicated situations as (i) and (ii) in this example.

\end{ex}

\end{document}